# An adaptive significance threshold criterion for massive multiple hypotheses testing

Cheng Cheng[1],*

*St. Jude Children's Research Hospital*

**Abstract:** This research deals with massive multiple hypothesis testing. First regarding multiple tests as an estimation problem under a proper population model, an error measurement called *Erroneous Rejection Ratio* (ERR) is introduced and related to the False Discovery Rate (FDR). ERR is an error measurement similar in spirit to FDR, and it greatly simplifies the analytical study of error properties of multiple test procedures. Next an improved estimator of the proportion of true null hypotheses and a data adaptive significance threshold criterion are developed. Some asymptotic error properties of the significant threshold criterion is established in terms of ERR under distributional assumptions widely satisfied in recent applications. A simulation study provides clear evidence that the proposed estimator of the proportion of true null hypotheses outperforms the existing estimators of this important parameter in massive multiple tests. Both analytical and simulation studies indicate that the proposed significance threshold criterion can provide a reasonable balance between the amounts of false positive and false negative errors, thereby complementing and extending the various FDR control procedures. S-plus/R code is available from the author upon request.

## 1. Introduction

The recent advancement of biological and information technologies made it possible to generate unprecedented large amounts of data for just a single study. For example, in a genome-wide investigation, expressions of tens of thousands genes and markers can be generated and surveyed simultaneously for their association with certain traits or biological conditions of interest. Statistical analysis in such applications poses a massive multiple hypothesis testing problem. The traditional approaches to controlling the probability of family-wise type-I error have proven to be too conservative in such applications. Recent attention has been focused on the control of *false discovery rate* (FDR) introduced by Benjamini and Hochberg [4]. Most of the recent methods can be broadly characterized into several approaches. Mixture-distribution partitioning [2, 24, 25] views the P values as random variables and models the P value distribution to generate estimates of the FDR levels at various significance levels. Significance analysis of microarrays (SAM; [32, 35]) employs permutation tests to inference simultaneously on order statistics. Empirical Baysian approaches include for example [10, 11, 17, 23, 28]. Tsai *et al.* [34] proposed models and estimators of the conditional FDR, and Bickel [6] takes a decision-theoretic approach. Recent theoretical developments on FDR control include Genovese and Wasserman [13, 14], Storey *et al.* [31], Finner and Roberts [12], and Abramovich *et al.* [1]. Recent theoretical development on control of generalized family-wise type-I error includes van der Laan *et al.* [36, 37], Dudoit *et al.* [9], and the references therein.

*Supported in part by the NIH grants U01 GM-061393 and the Cancer Center Support Grant P30 CA-21765, and the American Lebanese and Syrian Associated Charities (ALSAC).
[1]Department of Biostatistics, Mail Stop 768, St. Jude Childrens Research Hospital, 332 North Lauderdale Street, Memphis, TN 38105-2794, USA, e-mail: cheng.cheng@stjude.org

Benjamini and Hochberg [4] argue that as an alternative to the family-wise type-I error probability, FDR is a proper measurement of the amount of false positive errors, and it enjoys many desirable properties not possessed by other intuitive or heuristic measurements. Furthermore they develop a procedure to generate a significance threshold (P value cutoff) that guarantees the control of FDR under a pre-specified level. Similar to a significance test, FDR control requires one to specify a control level *a priori*. Storey [29] takes the point of view that in discovery-oriented applications neither the FDR control level nor the significance threshold may be specified before one sees the data (P values), and often the significance threshold is so determined *a posteriori* that allows for some "discoveries" (rejecting one or more null hypotheses). These "discoveries" are then scrutinized in confirmation and validation studies. Therefore it would be more appropriate to measure the false positive errors conditional on having rejected some null hypotheses, and for this purpose the positive FDR (pFDR; Storey [29]) is a meaningful measurement. Storey [29] introduces estimators of FDR and pFDR, and the concept of $q$-value which is essentially a neat representation of Benjamini and Hochberg's ([4]) step-up procedure possessing a Bayesian interpretation as the posterior probability of the null hypothesis ([30]). Reiner *et al.* [26] introduce the "FDR-adjusted P value" which is equivalent to the $q$-value. The $q$-value plot ([33]) allows for visualization of FDR (or pFDR) levels in relationship to significance thresholds or numbers of null hypotheses to reject. Other closely related procedures are the adaptive FDR control by Benjamini and Hochberg [3], and the recent two-stage linear step-up procedure by Benjamini *et al.* [5] which is shown to provide sure FDR control at any pre-specified level.

In discovery-oriented exploratory studies such as genome-wide gene expression survey or association rule mining in marketing applications, it is desirable to strike a meaningful balance between the amounts false positive and false negative errors than to control the FDR or pFDR alone. Cheng *et al.* [7] argue that it is not always clear in practice how to specify the threshold for either the FDR level or the significance level. Therefore, additional statistical guidelines beyond FDR control procedures are desirable. Genovese and Wasserman [13] extend FDR control to a minimization of the "false nondiscovery rate" (FNR) under a penalty of the FDR, i.e., $FNR + \lambda FDR$, where the penalty $\lambda$ is assumed to be specified *a priori*. Cheng *et al.* [7] propose to extract more information from the data (P values) and introduce three data-driven criteria for determination of the significance threshold.

This paper has two related goals: (1) develop a more accurate estimator of the proportion of true null hypotheses, which is an important parameter in all multiple hypothesis testing procedures; and (2) further develop the "profile information criterion" $I_p$ introduced in [7] by constructing a more data-adaptive criterion and study its asymptotic error behavior (as the number of tests tends to infinity) theoretically and via simulation. For theoretical and methodological development, a new meaningful measurement of the quantity of false positive errors, the *erroneous rejection ratio* (ERR), is introduced. Just like FDR, ERR is equal to the family-wise type-I error probability when all null hypotheses are true. Under the ergodicity conditions used in recent studies ([14, 31]), ERR is equal to FDR at any significant threshold (P value cut-off). On the other hand, ERR is much easier to handle analytically than FDR under distributional assumptions more widely satisfied in applications. Careful examination of each component in ERR gives insights into massive multiple testing in terms of the ensemble behavior of the P values. Quantities derived from ERR suggest to construct improved estimators of the null proportion (or the number of true null hypotheses) considered in [3, 29, 31], and the construction of an adaptive significance threshold criterion. The theoretical results demonstrate how the criterion can be calibrated with the Bonferroni adjustment to provide control of family-wise type-I error probability when all null hypotheses are true, and how the criterion behaves asymptotically, giving cautions and remedies in practice. The simulation results are consistent with the theory, and demonstrate that the proposed adaptive significance criterion is a useful and effective procedure complement to the popular FDR control methods.

This paper is organized as follows: Section 2 contains a brief review of FDR and the introduction of ERR; section 3 contains a brief review of the estimation of the proportion of null hypotheses, and the development of an improved estimator; section 4 develops the adaptive significance threshold criterion and studies its asymptotic error behavior (as the number of hypotheses tends to infinity) under proper distributional assumptions on the P values; section 5 contains a simulation study; and section 6 contains concluding remarks.

**Notation.** Henceforth, $\mathbb{R}$ denotes the real line; $\mathbb{R}^k$ denotes the $k$ dimensional Euclidean space. The symbol $\|\cdot\|_p$ denotes the $L^p$ or $\ell^p$ norm, and := indicates equal by definition. Convergence and convergence in probability are denoted by $\longrightarrow$ and $\longrightarrow_p$ respectively. A random variable is usually denoted by an upper-case letter such as $P$, $R$, $V$, etc. A cumulative distribution function (cdf) is usually denoted by $F$, $G$ or $H$; an empirical distribution function (EDF) is usually indicated by a tilde, e.g., $\widetilde{F}$. A population parameter is usually denoted by a lower-case Greek letter and a hat indicates an estimator of the parameter, e.g., $\widehat{\theta}$. Equivalence is denoted by $\simeq$, e.g., "$a_n \simeq b_n$ as $n \longrightarrow \infty$" means $\lim_{n \longrightarrow \infty} a_n/b_n = 1$.

## 2. False discovery rate and erroneous rejection ratio

Consider testing $m$ hypothesis pairs $(H_{0i}, H_{Ai})$, $i = 1, \ldots, m$. In many recent applications such as analysis of microarray gene differential expressions, $m$ is typically on the order of $10^5$. Suppose $m$ P values, $P_1, \ldots, P_m$, one for each hypothesis pair, are calculated, and a decision on whether to reject $H_{0i}$ is to be made. Let $m_0$ be the number of true null hypotheses, and let $m_1 := m - m_0$ be the number of true alternative hypotheses. The outcome of testing these $m$ hypotheses can be tabulated as in Table 1 (Benjamini and Hochberg [4]), where $V$ is the number of null hypotheses erroneously rejected, $S$ is the number of alternative hypotheses correctly captured, and $R$ is the total number of rejections.

Clearly only $m$ is known and only $R$ is observable. At least one *family-wise type-I error* is committed if $V > 0$, and procedures for multiple hypothesis testing have traditionally been produced for solely controlling the family-wise type-I error probability $\Pr(V > 0)$. It is well-known that such procedures often lack statistical power. In an effort to develop more powerful procedures, Benjamini and Hochberg ([4]) approached the multiple testing problem from a different perspective and introduced the concept of *false discovery rate* (FDR), which is, loosely speaking, the expected value of the ratio $V/R$. They introduced a simple and effective procedure for controlling the FDR under any pre-specified level.

It is convenient both conceptually and notationally to regard multiple hypotheses testing as an estimation problem ([7]). Define the parameter $\Theta = [\theta_1, \ldots, \theta_m]$ as $\theta_i = 1$ if $H_{Ai}$ is true, and $\theta_i = 0$ if $H_{0i}$ is true ($i = 1, \ldots, m$). The data consist of the P values $\{P_1, \ldots, P_m\}$, and under the assumption that each test is exact and unbiased, the population is described by the following probability model:

$$
\begin{aligned}
& P_i \sim \mathcal{P}_{i,\theta_i}; \\
& \mathcal{P}_{i,0} \text{ is } U(0,1), \text{ and } \mathcal{P}_{i,1} <_{st} U(0,1); \\
& \text{each } \mathcal{P}_{i,1} \text{ has a twice continuously differentiable cdf } F_i(\cdot),
\end{aligned}
\tag{2.1}
$$

for $i = 1, \ldots, m$, where $<_{st}$ stands for "stochastically less than." The P values are dependent in general and have a joint distribution on $\mathbb{R}^m$. By this model the marginal cdf of $P_i$ can be written as $G_i(t) = (1 - \theta_i)t + \theta_i F_i(t)$. Note $F_i(t) \geq t$ and $G_i(t) \geq t$ for $t \in [0,1]$.

A rejection procedure is an estimator of $\Theta$: $\widehat{\Theta} = \widehat{\Theta}(P_1, \ldots, P_m) = [\widehat{\theta}_1, \ldots, \widehat{\theta}_m] \in \{0,1\}^m$, where $\widehat{\theta}_i = 1$ indicates rejecting $H_{0i}$ in favor of $H_{Ai}$, $i = 1, \ldots, m$. With this notation, the random variables in Table 1 can be expressed as

$$
V = V_\Theta(\widehat{\Theta}) = \sum_{i=1}^m (1 - \theta_i)\widehat{\theta}_i;\ S = S_\Theta(\widehat{\Theta}) = \sum_{i=1}^m \theta_i \widehat{\theta}_i;\ R = R(\widehat{\Theta}) = \sum_{i=1}^m \widehat{\theta}_i.
\tag{2.2}
$$

A natural and perhaps the simplest procedure is the "hard-thresholding" (HT) estimator $\widehat{\Theta} = \widehat{\Theta}(\alpha)$ defined as

$$
HT(\alpha): \quad \widehat{\theta}_i = 1 \text{ iff } P_i \leq \alpha,
\tag{2.3}
$$

where $\alpha \in (0,1)$ is a significance threshold common to all hypotheses. Clearly for this procedure the random variables $V$, $S$, and $R$ all depend on $\alpha$. Traditional control of family-wise type-I error probability seeks to determine $\alpha$ so that $\Pr(V > 0) \leq \alpha^*$ for pre-specified $\alpha^*$. Genovese and Wasserman [14] list several procedures to determine $\alpha$. Benjamini and Hochberg [4] introduce a simple procedure to determine $\alpha$ so that the FDR is controlled at a given level.

## 2.1. False discovery rate and its control

The FDR as defined by Benjamini and Hochberg ([4]) can be expressed as

$$
FDR_\Theta(\widehat{\Theta}) = E\left[\frac{\sum_{i=1}^m \widehat{\theta}_i(1 - \theta_i)}{\sum_{i=1}^m \widehat{\theta}_i + \prod_{i=1}^m (1 - \widehat{\theta}_i)}\right],
\tag{2.4}
$$

which is equivalent to $E[V/R | R > 0]\Pr(R > 0)$. Let $P_{1:m} \leq P_{2:m} \leq \cdots \leq P_{m:m}$ be the order statistics of the P values, and let $\pi_0 = m_0/m$. Benjamini and Hochberg

TABLE 1
*Outcome tabulation of multiple hypotheses testing.*

| True Hypotheses | Rejected | Not Rejected | Total |
|---|---|---|---|
| $H_0$ | $V$ | $m_0 - V$ | $m_0$ |
| $H_A$ | $S$ | $m_1 - S$ | $m_1$ |
| Total | $R$ | $m - R$ | $m$ |

([4]) prove that for any specified $q^* \in (0,1)$ rejecting all the null hypotheses corresponding to $P_{1:m}, \ldots, P_{k^*:m}$ with $k^* = \max\{k : P_{k:m}/(k/m) \leq q^*\}$ controls the FDR at the level $\pi_0 q^*$, i.e., $FDR_\Theta(\widehat{\Theta}(P_{k^*:m})) \leq \pi_0 q^* \leq q^*$. Note this procedure is equivalent to applying the data-driven threshold $\alpha = P_{k^*:m}$ to all P values in (2.3), i.e., $HT(P_{k^*:m})$.

Recognizing the potential of constructing less conservative FDR controls by the above procedure, Benjamini and Hochberg ([3]) propose an estimator of $m_0$, $\widehat{m}_0$, (hence an estimator of $\pi_0$, $\widehat{\pi}_0 = \widehat{m}_0/m$), and replace $k/m$ by $k/\widehat{m}_0$ in determining $k^*$. They call this procedure "adaptive FDR control." The estimator $\widehat{\pi}_0 = \widehat{m}_0/m$ will be discussed in Section 3. A recent development in adaptive FDR control can be found in Benjamini *et al.* [5].

Similar to a significance test, the above procedure requires the specification of an FDR control level $q^*$ before the analysis is conducted. Storey ([29]) takes the point of view that for more discovery-oriented applications the FDR level is not specified *a priori*, but rather determined after one sees the data (P values), and it is often determined in a way allowing for some "discovery" (rejecting one or more null hypotheses). Hence a concept similar to, but different than FDR, the *positive false discovery rate* (pFDR) $E\big[V/R\big|R > 0\big]$, is more appropriate. Storey ([29]) introduces estimators of $\pi_0$, the FDR, and the pFDR from which the *q-values* are constructed for FDR control. Storey *et al.* ([31]) demonstrate certain desirable asymptotic conservativeness of the $q$-values under a set of ergodicity conditions.

### 2.2. Erroneous rejection ratio

As discussed in [3, 4], the FDR criterion has many desirable properties not possessed by other intuitive alternative criteria for multiple tests. In order to obtain an analytically convenient expression of FDR for more in-depth investigations and extensions, such as in [13, 14, 29, 31], certain fairly strong ergodicity conditions have to be assumed. These conditions make it possible to apply classical empirical process methods to the "FDR process." However, these conditions may be too strong for more recent applications, such as genome-wide tests for gene expression–phenotype association using microarrays, in which a substantial proportion of the tests can be strongly dependent. In such applications it may not be even reasonable to assume that the tests corresponding to the true null hypotheses are independent, an assumption often used in FDR research. Without these assumptions however, the FDR becomes difficult to handle analytically. An alternative error measurement in the same spirit of FDR but easier to handle analytically is defined below.

Define the *erroneous rejection ratio* (ERR) as

$$(2.5) \qquad ERR_\Theta(\widehat{\Theta}) = \frac{E[V_\Theta(\widehat{\Theta})]}{E[R(\widehat{\Theta})]} \Pr(R(\widehat{\Theta}) > 0).$$

Just like FDR, when all null hypotheses are true $ERR = \Pr(R(\widehat{\Theta}) > 0)$, which is the family-wise type-I error probability because now $V_\Theta(\widehat{\Theta}) = R(\widehat{\Theta})$ with probability one. Denote by $V(\alpha)$ and $R(\alpha)$ respectively the $V$ and $R$ random variables and by $ERR(\alpha)$ the ERR for the hard-thresholding procedure $HT(\alpha)$; thus

$$(2.6) \qquad ERR(\alpha) = \frac{E[V(\alpha)]}{E[R(\alpha)]} \Pr(R(\alpha) > 0).$$

Careful examination of each component in $ERR(\alpha)$ reveals insights into multiple tests in terms of the ensemble behavior of the P values. Note

$$E[V(\alpha)] = \sum_{i=1}^{m}(1-\theta_i)\Pr(\widehat{\theta}_i = 1) = m_0\alpha$$
$$E[R(\alpha)] = \sum_{i=1}^{m}\Pr(\widehat{\theta}_i = 1) = m_0\alpha + \sum_{j:\theta_j=1} F_j(\alpha)$$
$$\Pr(R(\alpha) > 0) = \Pr(P_{1:m} \leq \alpha).$$

Define $H_m(t) := m_1^{-1}\sum_{j:\theta_j=1} F_j(t)$ and $F_m(t) := m^{-1}\sum_{i=1}^{m} G_i(t) = \pi_0 t + (1-\pi_0)H_m(t)$. Then

$$(2.7) \qquad ERR(\alpha) = \frac{\pi_0 \alpha}{F_m(\alpha)}\Pr(P_{1:m} \leq \alpha).$$

The functions $H_m(\cdot)$ and $F_m(\cdot)$ both are cdf's on $[0,1]$; $H_m$ is the average of the P value marginal cdf's corresponding to the true alternative hypotheses, and $F_m$ is the average of all P value marginal cdf's. $F_m$ describes the ensemble behavior of all P values and $H_m$ describes the ensemble behavior of the P values corresponding to the true alternative hypotheses. Cheng et al. ([7]) observe that the EDF of the P values $\widetilde{F}_m(t) := m^{-1}\sum_{i=1}^{m} I(P_i \leq t)$, $t \in \mathbb{R}$ is an unbiased estimator of $F_m(\cdot)$, and if the tests $\widehat{\theta}$ ($i = 1, \ldots, m$) are not strongly correlated asymptotically in the sense that $\sum_{i \neq j} Cov(\widehat{\theta}_i, \widehat{\theta}_j) = o(m^2)$ as $m \longrightarrow \infty$, $\widetilde{F}_m(\cdot)$ is "asymptotically consistent" for $F_m$ in the sense that $|\widetilde{F}_m(t) - F_m(t)| \longrightarrow_p 0$ for every $t \in \mathbb{R}$. This prompts possibilities for the estimation of $\pi_0$, data-adaptive determination of $\alpha$ for the $HT(\alpha)$ procedure, and the estimation of FDR. The first two will be developed in detail in subsequent sections. Cheng et al. ([7]) and Pounds and Cheng ([25]) develop smooth FDR estimators.

Let $FDR(\alpha) := E[V(\alpha)/R(\alpha)|R(\alpha) > 0]\Pr(R(\alpha) > 0)$. $ERR(\alpha)$ is essentially $FDR(\alpha)$. Under the hierarchical (or random effect) model employed in several papers ([11, 14, 29, 31]), the two quantities are equivalent, that is, $FDR(\alpha) = ERR(\alpha)$ for all $\alpha \in (0,1]$, following from Lemma 2.1 in [14]. More generally $ERR/FDR = \{E[V]/E[R]\}/E[V/R|R > 0]$ provided $\Pr(R > 0) > 0$. Asymptotically as $m \longrightarrow \infty$, if $\Pr(R > 0) \longrightarrow 1$ then $E[V/R|R > 0] \simeq E[V/R]$; if furthermore $E[V/R] \simeq E[V]/E[R]$, then $ERR/FDR \longrightarrow 1$. Identifying reasonable sufficient (and necessary) conditions for $E[V/R] \simeq E[V]/E[R]$ to hold remains an open problem at this point.

Analogous to the relationship between FDR and pFDR, define the *positive ERR*, $pERR := E[V]/E[R]$. Both quantities are well-defined provided $\Pr(R > 0) > 0$. The relationship between pERR and pFDR is the same as that between ERR and FDR described above.

The error behavior of a given multiple test procedure can be investigated in terms of either FDR (pFDR) or ERR (pERR). The ratio $pERR = E[V]/E[R]$ can be handled easily under arbitrary dependence among the tests because $E[V]$ and $E[R]$ are simply means of sums of indicator random variables. The only possible challenging component in $ERR(\alpha)$ is $\Pr(R(\alpha) > 0) = \Pr(P_{1:m} \leq \alpha)$; some assumptions on the dependence among the tests has to be made to obtain a concrete analytical form for this probability, or an upper bound for it. Thus, as demonstrated in Section 4, ERR is an error measurement that is easier to handle than FDR under more complex and application-pertinent dependence among the tests, in assessing analytically the error properties of a multiple hypothesis testing procedure.

A fine technical point is that FDR (pFDR) is always well-defined and ERR (pERR) is always well-defined under the convention $a \cdot 0 = 0$ for $a \in [-\infty, +\infty]$.

Compared to FDR (pFDR), ERR (pERR) is slightly less intuitive in interpretation. For example, FDR can be interpreted as the expected proportion of false positives among all positive findings, whereas ERR can be interpreted as the proportion of the number of false positives expected out of the total number of positive findings expected. Nonetheless, ERR (pERR) is still of practical value given its close relationship to FDR (pFDR), and is more convenient to use in analytical assessments of a multiple test procedure.

## 3. Estimation of the proportion of null hypotheses

The proportion of the true null hypotheses $\pi_0$ is an important parameter in all multiple test procedures. A delicate component in the control or estimation of FDR (or ERR) is the estimation of $\pi_0$. The cdf $F_m(t) = \pi_0 t + (1 - \pi_0) H_m(t)$, $t \in [0, 1]$, along with the fact that the EDF $\widetilde{F}_m$ is its unbiased estimator provides a clue for estimating $\pi_0$. Because for any $t \in (0, 1)$ $\pi_0 = [H_m(t) - F_m(t)]/[H_m(t) - t]$, a plausible estimator of $\pi_0$ is

$$\widehat{\pi}_0 = \frac{\Lambda - \widetilde{F}_m(t_0)}{\Lambda - t_0}$$

for properly chosen $\Lambda$ and $t_0$. Let $Q_m(u) := F_m^{-1}(u)$, $u \in [0, 1]$ be the quantile function of $F_m$ and let $\widetilde{Q}_m(u) := \widetilde{F}_m^{-1}(u) := \inf\{x : \widetilde{F}_m(x) \geq u\}$ be the *empirical quantile function* (EQF), then $\pi_0 = [H_m(Q_m(u)) - u]/[H_m(Q_m(u)) - Q_m(u)]$, for $u \in (0, 1)$, and with $\Lambda_1$ and $u_0$ properly chosen

$$\widehat{\pi}_0 = \frac{\Lambda_1 - u_0}{\Lambda_1 - \widetilde{Q}_m(u_0)}$$

is a plausible estimator. The existing $\pi_0$ estimators take either of the above representations with minor modifications.

Clearly it is necessary to have $\Lambda_1 \geq u_0$ for a meaningful estimator. Because $Q_m(u_0) \leq u_0$ by the stochastic order assumption [cf. (2.1)], choosing $\Lambda_1$ too close to $u_0$ will produce an estimator much biased downward. Benjamini and Hochberg ([3]) use the heuristic that if $u_0$ is so chosen that all P values corresponding to the alternative hypotheses concentrate in $[0, Q_m(u_0)]$ then $H_m(Q_m(u_0)) = 1$; thus setting $\Lambda_1 = 1$. Storey ([29]) uses a similar heuristic to set $\Lambda = 1$.

### 3.1. Existing estimators

Taking a graphical approach Schweder and Spjøtvoll [27] propose an estimator of $m_0$ as $\widehat{m}_0 = m(1 - \widetilde{F}_m(\lambda))/(1 - \lambda)$ for a properly chosen $\lambda$; hence a corresponding estimator of $\pi_0$ is $\widehat{\pi}_0 = \widehat{m}_0/m = (1 - \widetilde{F}_m(\lambda))/(1 - \lambda)$. This is exactly Storey's ([29]) estimator. Storey observes that $\lambda$ is a tuning parameter that dictates the bias and variance of the estimator, and proposes computing $\widehat{\pi}_0$ on a grid of $\lambda$ values, smoothing them by a spline function, and taking the smoothed $\widehat{\pi}_0$ at $\lambda = 0.95$ as the final estimator. Storey *et al.* ([31]) propose a bootstrap procedure to estimate the mean-squared error (MSE) and pick the $\lambda$ that gives the minimal estimated MSE. It will be seen in the simulation study (Section 5) that this estimator tends to be biased downward.

Approaching to the problem from the quantile perspective Benjamini and Hochberg ([3]) propose $\widehat{m}_0 = \min\{1 + (m+1-j)/(1-P_{j:m}), m\}$ for a properly chosen $j$; hence

$$\widehat{\pi}_0 = \min\left\{\frac{1}{m} + \left[\frac{1-P_{j:m}}{1-j/m+1/m}\right]^{-1}, 1\right\}.$$

The index $j$ is determined by examining the slopes $S_i = (1 - P_{i:m})/(m+1-i)$, $i = 1, \ldots, m$, and is taken to be the smallest index such that $S_j < S_{j-1}$. Then $\widehat{m}_0 = \min\{1 + 1/S_j, m\}$. It is not difficult to see why this estimator tends to be too conservative (i.e., too much biased upward): as $m$ gets large the event $\{S_j < S_{j-1}\}$ tends to occur early (i.e., at small $j$) with high probability. By definition, $S_j < S_{j-1}$ if and only if

$$\frac{1-P_{j:m}}{m+1-j} < \frac{1-P_{j-1:m}}{m+2-j},$$

if and only if

$$P_{j:m} > \frac{1}{m+2-j} + \frac{m+1-j}{m+2-j}P_{j-1:m}.$$

Thus, as $m \to \infty$,

$$\Pr(S_j < S_{j-1}) = \Pr\left(P_{j:m} > \frac{1}{m+2-j} + \frac{m+1-j}{m+2-j}P_{j-1:m}\right) \longrightarrow 1,$$

for fixed or small enough $j$ satisfying $j/m \longrightarrow \delta \in [0,1)$. The conservativeness will be further demonstrated by the simulation study in Section 5.

Recently Mosig *et al.* ([21]) proposed an estimator of $m_0$ by a recursive algorithm, which is clarified and shown by Nettleton and Hwang [22] to converge under a fixed partition (histogram bins) of the P value order statistics. In essence the algorithm searches in the right tail of the P value histogram to determine a "bend point" when the histogram begins to become flat, and then takes this point for $\lambda$ (or $j$).

For a two-stage adaptive control procedure Benjamini *et al.* ([5]) consider an estimator of $m_0$ derived from the first-stage FDR control at the more conservative $q/(1+q)$ level than the targeted control level $q$. Their simulation study indicates that with comparable bias this estimator is much less variable than the estimators by Benjamini and Hochberg [3] and Storey *et al.* [31], thus possessing better accuracy. Recently Langaas *et al.* ([19]) proposed an estimator based on nonparametric estimation of the P value density function under monotone and convex contraints.

### 3.2. An estimator by quantile modeling

Intuitively, the stochastic order requirement in the distributional model (2.1) implies that the cdf $F_m(\cdot)$ is approximately concave and hence the quantile function $Q_m(\cdot)$ is approximately convex. When there is a substantial proportion of true null and true alternative hypotheses, there is a "bend point" $\tau_m \in (0,1)$ such that $Q_m(\cdot)$ assumes roughly a nonlinear shape in $[0, \tau_m]$, primarily dictated by the distributions of the P values corresponding to the true alternative hypotheses, and $Q_m(\cdot)$ is essentially linear in $[\tau_m, 1]$, dictated by the $U(0,1)$ distribution for the null P values. The estimation of $\pi_0$ can benefit from properly capturing this shape characteristic by a model.

Clearly $\pi_0 \leq [1 - \tau_m]/[H_m(Q_m(\tau_m)) - Q_m(\tau_m)]$. Again heuristically if all P values corresponding to the alternative hypotheses concentrate in $[0, Q_m(\tau_m)]$, then $H_m(Q_m(\tau_m)) = 1$. A strategy then is to construct an estimator of $Q_m(\cdot)$, $\widehat{Q}_m^*(\cdot)$, that possesses the desirable shape described above, along with a bend point $\widehat{\tau}_m$, and set

$$\widehat{\pi}_0 = \frac{1 - \widehat{\tau}_m}{1 - \widehat{Q}_m^*(\widehat{\tau}_m)}, \tag{3.1}$$

which is the inverse slope between the points $(\widehat{\tau}_m, \widehat{Q}_m^*(\widehat{\tau}_m))$ and $(1,1)$ on the unit square.

Model (2.1) implies that $Q_m(\cdot)$ is twice continuously differentiable. Taylor expansion at $t=0$ gives $Q_m(t) = q_m(0)t + \frac{1}{2}q_m'(\xi_t)t^2$ for $t$ close to 0 and $0 < \xi_t < t$, where $q_m(\cdot)$ is the first derivative of $Q_m(\cdot)$, i.e., the *quantile density function* (qdf), and $q_m'(\cdot)$ is the second derivative of $Q_m(\cdot)$. This suggests the following definition (model) of an approximation of $Q_m$ by a convex, two-piece function joint smoothly at $\tau_m$. Define $\underline{Q}_m(t) := \min\{Q_m(t), t\}$, $t \in [0,1]$, define the bend point $\tau_m := \operatorname{argmax}_t\{t - \underline{Q}_m(t)\}$ and assume that it exists uniquely, with the convention that $\tau_m = 0$ if $Q_m(t) = t$ for all $t \in [0,1]$. Define

$$Q_m^*(t; \gamma, a, d, b_1, b_0, \tau_m) = \begin{cases} at^\gamma + dt, & 0 \leq t \leq \tau_m \\ b_0 + b_1 t, & t \geq \tau_m \end{cases} \tag{3.2}$$

where

$$b_1 = [1 - Q_m(\tau_m)] / (1 - \tau_m)$$
$$b_0 = 1 - b_1 = [Q_m(\tau_m) - \tau_m] / (1 - \tau_m),$$

and $\gamma$, $a$ and $d$ are determined by minimizing $\|Q_m^*(\cdot; \gamma, a, d, b_1, b_0, \tau_m) - Q_m(\cdot)\|_1$ under the following constraints:

$$\begin{cases} \gamma \geq 1, \ a \geq 0, \ 0 \leq d \leq 1 \\ \gamma = a = 1, \ d = 0 \text{ if and only if } \tau_m = 0 \\ a\tau_m^\gamma + d\tau_m = b_0 + b_1\tau_m \ (\text{continuity at } \tau_m) \\ a\gamma\tau_m^{\gamma-1} + d = b_1 \ (\text{smoothness at } \tau_m). \end{cases}$$

These constraints guarantee that the two pieces are joint smoothly at $\tau_m$ to produce a convex and continuously differentiable quantile function that is the closest to $Q_m$ on $[0,1]$ in the $L^1$ norm, and that there is no over-parameterization if $Q_m$ coincides with the 45-degree line. $Q_m^*$ will be called the *convex backbone* of $Q_m$.

The smoothness constraints force $a$, $d$ and $\gamma$ to be interdependent via $b_0$, $b_1$ and $\tau_m$. For example,

$$\begin{cases} a = a(\gamma) = -b_0 / [(\gamma - 1)\tau_m] \ (\text{for } \gamma > 1) \\ d = d(\gamma) = b_1 - a(\gamma)\gamma\tau_m^{\gamma-1}. \end{cases}$$

Thus the above constrained minimization is equivalent to

$$\min_\gamma \|Q_m^*(\cdot; \gamma, a(\gamma), d(\gamma), b_1, b_0, \tau_m) - Q_m(\cdot)\|_1 \tag{3.3}$$

$$\text{subject to } \begin{cases} \gamma \geq 1, \ a(\gamma) \geq 0, \ 0 \leq d(\gamma) \leq 1 \\ \gamma = a = 1, \ d = 0 \text{ if and only if } \tau_m = 0. \end{cases}$$

An estimator of $\pi_0$ is obtained by plugging an estimator of the convex backbone $Q_m^*$, $\widehat{Q}_m^*$, into (3.1). The convex backbone can be estimated by replacing $Q_m$ with the EQF $\widetilde{Q}_m$ in the above process. However, instead of using the raw EQF, the estimation can benefit from properly smoothing the modified EQF $\underline{\widetilde{Q}}_m(t) := \min\{\widetilde{Q}_m(t), t\}$, $t \in [0,1]$ into a smooth and approximately convex EQF, $\widehat{Q}_m(\cdot)$. This smooth and approximately convex EQF can be obtained by repeatedly smoothing the modified EQF $\underline{\widetilde{Q}}_m(\cdot)$ by the variation-diminishing spline (VD-spline; de Boor [7], P.160). Denote by $B_{j,\mathbf{t},k}$ the $j$th order-$k$ B spline with extended knot sequence $\mathbf{t} = t_1, \ldots, t_{n+k}$ ($t_1 = \ldots t_k = 0 < t_{k+1} < \ldots < t_n < t_{n+1} = \ldots = t_{n+k} = 1$) and $t_j^* := \sum_{\ell=j+1}^{j+k-1} t_\ell/(k-1)$. The VD-spline approximation of a function $h: [0,1] \to \mathbb{R}$ is defined as

$$(3.4) \qquad \widehat{h}(u) := \sum_{j=1}^{n} h(t_j^*) B_{j,\mathbf{t};k}(u), \qquad u \in [0,1].$$

The current implementation takes $k = 5$ (thus quartic spline for $\widehat{Q}_m$ and cubic spline for its derivative, $\widehat{q}_m$), and sets the interior knots in $\mathbf{t}$ to the ordered unique numbers in $\{\frac{1}{m}, \frac{2}{m}, \frac{3}{m}, \frac{4}{m}\} \cup \{\widetilde{F}_m(t), \ t = 0.001, 0.003, 0.00625, 0.01, 0.0125, 0.025, 0.05, 0.1, 0.25\}$. The knot sequence is so designed that the variation in the quantile function in a neighborhood close to zero (corresponding to small P values) can be well captured; whereas the right tail (corresponding to large P values) is greatly smoothed. Key elements in the algorithm, such as the interior knots positions, the $t_j^*$ positions, etc., are illustrated in Figure 1.

Upon obtaining the smooth and approximately convex EQF $\widehat{Q}_m(\cdot)$, the convex backbone estimator $\widehat{Q}_m^*(\cdot)$ is constructed by replacing $Q_m(\cdot)$ with $\widehat{Q}_m(\cdot)$ in (3.3) and numerically solving the optimization with a proper search algorithm. This algorithm produces the estimator $\widehat{\pi}_0$ in (3.1) at the same time.

Note that in general the parameters $\gamma$, $a$, $d$, $b_0$, $b_1$, $\pi_0 := m_0/m$, and their corresponding estimators all depend on $m$. For the sake of notational simplicity this dependency has been and continues to be suppressed in the notation. Furthermore, it is assumed that $\lim_{m\to\infty} m_0/m$ exists. For studying asymptotic properties, henceforth let $\{P_1, P_2, \ldots\}$ be an infinite sequence of P values, and let $\mathbf{P}_m := \{P_1, \ldots, P_m\}$.

## 4. Adaptive profile information criterion

### 4.1. The adaptive profile information criterion

We now develop an adaptive procedure to determine a significance threshold for the $HT(\alpha)$ procedure. The estimation perspective allows one to draw an analogy between multiple hypothesis testing and the classical variable selection problem: setting $\widehat{\theta}_i = 1$ (i.e., rejecting the $i$th null hypothesis) corresponds to including the $i$th variable in the model. A traditional model selection criterion such AIC usually consists of two terms, a model-fitting term and a penalty term. The penalty term is usually some measure of model complexity reflected by the number of parameters to be estimated. In the context of massive multiple testing a natural penalty (complexity) measurement would be the expected number of false positives $E[V(\alpha)] = \pi_0 m\alpha$ under model (2.1). When a parametric model is fully specified, the model-fitting term is usually a likelihood function or some similar quantity. In the context of massive multiple testing the stochastic order assumption in model (2.1) suggests using a proper quantity measuring the lack-off-fit from $U(0,1)$ in the ensemble distribution of the P values on the interval $[0, \alpha]$. Cheng et al. ([7]) considered such a measurement that is an $L^2$ distance. The concept of convex backbone facilitates the derivation of a measurement more adaptive to the ensemble distribution of the P values. Given the convex backbone $Q_m^*(\cdot) := Q_m^*(\cdot; \gamma, a, d, b_1, b_0, \tau_m)$ as defined in (3.2), the "model-fitting" term can be defined as the $L^\gamma$ distance between $Q_m^*(\cdot)$ and uniformity on $[0, \alpha]$:

$$D_\gamma(\alpha) := \left[ \int_0^\alpha (t - Q_m^*(t))^\gamma \, dt \right]^{1/\gamma}, \qquad \alpha \in (0, 1].$$

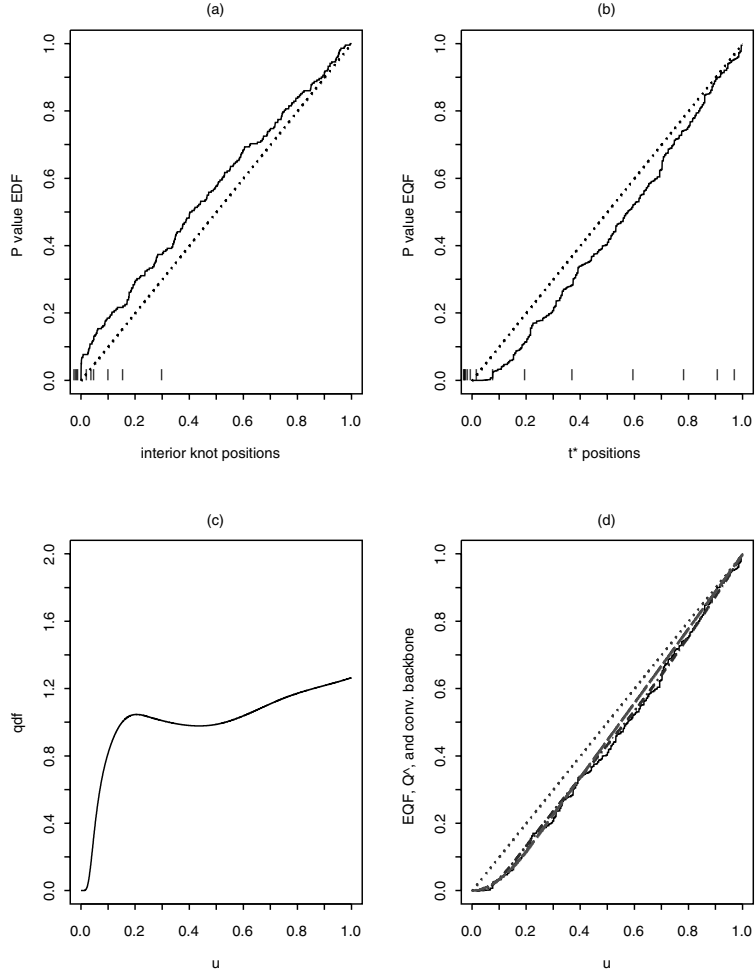

FIG 1. (a) *The interior knot positions indicated by | and the P value EDF;* (b) *the positions of $t_j^*$ indicated by | and the P value EQF;* (c) $\widehat{q}_m$: *the derivative of $\widehat{Q}_m$;* (d) *the P value EQF (solid), the smoothed EQF $\widehat{Q}_m$ from Algorithm 1 (dash-dot), and the convex backbone $\widehat{Q}_m^*$ (long dash).*

The adaptivity is reflected by the use of the $L^\gamma$ distance: Recall that the larger the $\gamma$, the higher concentration of small P values, and the norm inequality (Hardy, Littlewood, and Pólya [16], P.157) implies that $D_{\gamma_2}(\alpha) \geq D_{\gamma_1}(\alpha)$ for every $\alpha \in (0, 1]$ if $\gamma_2 > \gamma_1$.

Clearly $D_\gamma(\alpha)$ is non-decreasing in $\alpha$. Intuitively one possibility would be to maximize a criterion like $D_\gamma(\alpha) - \lambda \pi_0 m \alpha$. However, the two terms are not on the same order of magnitude when $m$ is very large. The problem is circumvented by using $1/D_\gamma(\alpha)$, which also makes it possible to obtain a closed-form solution to approximately optimizing the criterion.

Thus define the *Adaptive Profile Information* (API) criterion as

$$(4.1) \qquad API(\alpha) := \left[ \int_0^\alpha (t - Q_m^*(t))^\gamma \, dt \right]^{-1/\gamma} + \lambda(m, \pi_0, d) m \pi_0 \alpha,$$

for $\alpha \in (0,1)$ and $Q_m^*(\cdot) := Q_m^*(\cdot; \gamma, a, d, b_1, b_0, \tau_m)$ as defined in (3.2). One seeks to minimize $API(\alpha)$ to obtain an adaptive significance threshold for the $HT(\alpha)$ procedure.

With $\gamma > 1$, the integral can be approximated by $\int_0^\alpha ((1-d)t)^\gamma dt = (1-d)(\gamma+1)^{-1}\alpha^{\gamma+1}$. Thus

$$API(\alpha) \approx \overline{API}(\alpha) := (1-d)^{-1} \left[ \frac{1}{\gamma+1} \alpha^{\gamma+1} \right]^{-1/\gamma} + \lambda(m, \pi_0, d) m \pi_0 \alpha.$$

Taking the derivative of $\overline{API}(\cdot)$ and setting it to zero gives

$$\alpha^{-(2\gamma+1)/\gamma} = (1-d)(\gamma+1)^{-1/\gamma} \frac{\gamma}{\gamma+1} \lambda(m, \pi_0, d) m \pi_0.$$

Solving for $\alpha$ gives

$$\alpha^* = \left[ \frac{(\gamma+1)^{(1+1/\gamma)}}{(1-d)\pi_0 \gamma} \right]^{\gamma/(2\gamma+1)} [\lambda(m, \pi_0, d) m]^{-\gamma/(2\gamma+1)},$$

which is an approximate minimizer of $API$. Setting $\lambda(m, \pi_0, d) = m^{\beta \pi_0}/(1-d)$ and $\beta = 2\pi_0/\gamma$ gives

$$\alpha^* = \left[ \frac{(\gamma+1)^{(1+1/\gamma)}}{\pi_0 \gamma} \right]^{\gamma/(2\gamma+1)} m^{-(1+2\pi_0^2/\gamma)\gamma/(2\gamma+1)}.$$

This particular choice for $\lambda$ is motivated by two facts. When most of the P values have the $U(0,1)$ distribution (equivalently, $\pi_0 \approx 1$), the $d$ parameter of the convex backbone can be close to 1; thus with $1-d$ in the denominator, $\alpha^*$ can be unreasonably high in such a case. This issue is circumvented by putting $1-d$ in the denominator of $\lambda$, which eliminates $1-d$ from the denominator of $\alpha^*$. Next, it is instructive to compare $\alpha^*$ with the Bonferroni adjustment $\alpha_{Bonf}^* = \alpha_0/m$ for a pre-specified $\alpha_0$. If $\gamma$ is large, then $\alpha_{Bonf}^* < \alpha^* \approx O(m^{-1/2})$ as $m \longrightarrow \infty$. Although the derivation required $\gamma > 1$, $\alpha^*$ is still well defined even if $\pi_0 = 1$ (implying $\gamma = 1$), and in this case $\alpha^* = 4^{1/3} m^{-1}$ is comparable to $\alpha_{Bonf}^*$ as $m \longrightarrow \infty$. This in fact suggests the following significance threshold calibrated with the Bonferroni adjustment:

$$(4.2) \qquad \alpha_{cal}^* := 4^{-1/3} \left( \frac{\gamma}{\pi_0} \right) \alpha_0 \alpha^* = A(\pi_0, \gamma) m^{-B(\pi_0, \gamma)},$$

which coincides with the Bonferroni threshold $\alpha_0 m^{-1}$ when $\pi_0 = 1$, where

$$A(x, y) := \left[ y/(4^{1/3} x) \right] \alpha_0 \left[ (y+1)^{(1+1/y)}/(xy) \right]^{y/(2y+1)}$$

$$(4.3)$$

$$B(x, y) := (1 + 2x^2/y) y/(2y+1).$$

The factor $\alpha_0$ serves asymptotically as a calibrator of the adaptive significance threshold to the Bonferroni threshold in the least favorable scenario $\pi_0 = 1$, i.e., all null hypotheses are true. Analysis of the asymptotic ERR of the $HT(\alpha^*_{cal})$ procedure suggests a few choices of $\alpha_0$ in practice.

### 4.2. Asymptotic ERR of $HT(\alpha^*_{cal})$

Recall from (2.7) that

$$ERR(\alpha) = \left[\pi_0 \alpha / F_m(\alpha)\right] \Pr(P_{1:m} \leq \alpha).$$

The probability $\Pr(P_{1:m} \leq \alpha)$ is not tractable in general, but an upper bound can be obtained under a reasonable assumption on the set $\mathbf{P}_m$ of the $m$ P values. Massive multiple tests are mostly applied in exploratory studies to produce "inference-guided discoveries" that are either subject to further confirmation and validation, or helpful for developing new research hypotheses. For this reason often all the alternative hypotheses are two-sided, and hence so are the tests. It is instructive to first consider the case of $m$ two-sample $t$ tests. Conceptually the data consist of $n_1$ i.i.d. observations on $\mathbb{R}^m$ $\mathbf{X}_i = [X_{i1}, X_{i2}, \ldots, X_{im}]$, $i = 1, \ldots, n_1$ in the first group, and $n_2$ i.i.d. observations $\mathbf{Y}_i = [Y_{i1}, Y_{i2}, \ldots, Y_{im}]$, $i = 1, \ldots, n_2$ in the second group. The hypothesis pair $(H_{0k}, H_{Ak})$ is tested by the two-sided two-sample t statistic $T_k = |T(\mathcal{X}_k, \mathcal{Y}_k, n_1, n_2)|$ based on the data $\mathcal{X}_k = \{X_{1k}, \ldots, X_{n_1 k}\}$ and $\mathcal{Y}_k = \{Y_{1k}, \ldots, Y_{n_2 k}\}$. Often in biological applications that study gene signaling pathways (see e.g., Kuo *et al.* [18], and the simulation model in Section 5), $X_{ik}$ and $X_{ik'}$ ($i = 1, \ldots, n_1$) are either positively or negatively correlated for certain $k \neq k'$, and the same holds for $Y_{ik}$ and $Y_{ik'}$ ($i = 1, \ldots, n_2$). Such dependence in data raises positive association between the two-sided test statistics $T_k$ and $T_{k'}$ so that $\Pr(T_k \leq t | T'_k \leq t) \geq \Pr(T_k \leq t)$, implying $\Pr(T_k \leq t, T_{k'} \leq t) \geq \Pr(T_k \leq t) \Pr(T_{k'} \leq t)$, $t \geq 0$. Then the P values in turn satisfy $\Pr(P_k > \alpha, P_{k'} > \alpha) \geq \Pr(P_k > \alpha) \Pr(P_{k'} > \alpha)$, $\alpha \in [0, 1]$. It is straightforward to generalize this type of dependency to more than two tests. Alternatively, a direct model for the P values can be constructed.

**Example 4.1.** Let $\mathcal{J} \subseteq \{1, \ldots, m\}$ be a nonempty set of indices. Assume $P_j = P_0^{X_j}$, $j \in \mathcal{J}$, where $P_0$ follows a distribution $F_0$ on $[0, 1]$, and $X_j$'s are i.i.d. continuous random variables following a distribution $H$ on $[0, \infty)$, and are independent of the P values. Assume that the $P_i$'s for $i \notin \mathcal{J}$ are either independent or related to each other in the same fashion. This model mimics the effect of an activated gene signaling pathway that results in gene differential expression as reflected by the P values: the set $\mathcal{J}$ represents the genes involved in the pathway, $P_0$ represents the underlying activation mechanism, and $X_j$ represents the noisy response of gene $j$ resulting in $P_j$. Because $P_i > \alpha$ if and only if $X_j < \log \alpha / \log P_0$, direct calculations using independence of the $X_j$'s show that

$$\Pr\left(\bigcap_{j \in \mathcal{J}} \{P_j > \alpha\}\right) = \int_0^1 \Pr\left(\bigcap_{j \in \mathcal{J}} \left\{X_j < \frac{\log \alpha}{\log t}\right\}\right) dF_0(t) = E\left[\left[H\left(\frac{\log \alpha}{\log P_0}\right)\right]^{|\mathcal{J}|}\right],$$

where $|\mathcal{J}|$ is the cardinality $\mathcal{J}$. Next

$$\prod_{j \in \mathcal{J}} \Pr(P_j > \alpha) = \prod_{j \in \mathcal{J}} \int_0^1 \left[H\left(\frac{\log \alpha}{\log t}\right)\right] dF_0(t) = \left[E\left[H\left(\frac{\log \alpha}{\log P_0}\right)\right]\right]^{|\mathcal{J}|}.$$

Finally $\Pr\left(\cap_{j \in \mathcal{J}}\{P_j > \alpha\}\right) \geq \prod_{j \in \mathcal{J}} \Pr(P_j > \alpha)$, following from Jensen's inequality.

The above considerations lead to the following definition.

**Definition 4.1.** The set of P values $\mathbf{P}_m$ has the *positive orthant dependence property* if for any $\alpha \in [0, 1]$

$$\Pr\left(\bigcap_{i=1}^{m}\{P_i > \alpha\}\right) \geq \prod_{i=1}^{m} \Pr(P_i > \alpha).$$

This type of dependence is similar to the positive quadrant dependence introduced by Lehmann [20].

Now define the upper envelope of the cdf's of the P values as

$$\overline{F}_m(t) := \max_{i=1,\ldots,m}\{G_i(t)\}, \qquad t \in [0, 1],$$

where $G_i$ is the cdf of $P_i$. If $\mathbf{P}_m$ has the positive orthant dependence property then

$$\Pr(P_{1:m} \leq \alpha) = 1 - \Pr\left(\bigcap_{i=1}^{m}\{P_i > \alpha\}\right) \leq 1 - \prod_{i=1}^{m}\Pr(P_i > \alpha) \leq 1 - (1 - \overline{F}_m(\alpha))^m,$$

implying

(4.4) $$ERR(\alpha^*_{cal}) \leq \frac{\pi_0 \alpha^*_{cal}}{\pi_0 \alpha^*_{cal} + (-\pi_0) H_m(\alpha^*_{cal})} \left[1 - (1 - \overline{F}_m(\alpha^*_{cal}))^m\right].$$

Because $\alpha^*_{cal} \longrightarrow 0$ as $m \longrightarrow \infty$, the asymptotic magnitude of the above ERR can be established by considering the magnitude of $\overline{F}_m(t_m)$ and $H_m(t_m)$ as $t_m \longrightarrow 0$. The following definition makes this idea rigorous.

**Definition 4.2.** The set of $m$ P values $\mathbf{P}_m$ is said to be *asymptotically stable* as $m \longrightarrow \infty$ if there exists sequences $\{\beta_m\}$, $\{\eta_m\}$, $\{\psi_m\}$, $\{\xi_m\}$ and constants $\beta^*$, $\beta_*$, $\eta$, $\psi^*$, $\psi_*$, and $\xi$ such that

$$\overline{F}_m(t) \simeq \beta_m t^{\eta_m}, \quad H_m(t) \simeq \psi_m t^{\xi_m}, \quad t \longrightarrow 0$$

and

$$0 < \beta_* \leq \beta_m \leq \beta^* < \infty, \quad 0 < \eta \leq \eta_m \leq 1$$
$$0 < \psi_* \leq \psi_m \leq \psi^* < \infty, \quad 0 < \xi \leq \xi_m \leq 1$$

for sufficiently large $m$.

This definition essentially says that $\mathbf{P}_m$ is regarded as asymptotically stable if the ensemble distribution functions $F_m(\cdot)$ and $H_m(\cdot)$ vary in the left tail similarly to Beta distributions.

The following theorem establishes the asymptotic magnitude of an upper bound of $ERR(\alpha^*_{cal})$ – the ERR of applying $\alpha^*_{cal}$ in the hard-thresholding procedure (2.3).

**Theorem 4.1.** *Let $\psi_*$, $\xi_m$, $\beta^*$, and $\eta$ be as given in Definition 4.2, and let $A(\cdot, \cdot)$ and $B(\cdot, \cdot)$ be as defined in (4.3). If the set of P values $\mathbf{P}_m$ is asymptotically stable and has the positive orthant dependence property for sufficiently large $m$, then $ERR(\alpha^*_{cal}) \leq \Psi(\alpha^*_{cal})$, and $\Psi(\alpha^*_{cal})$ satisfies*

(a) *if $\pi_0 = 1$ for all $m$, $\lim_{m \to \infty} \Psi(\alpha^*_{cal}) = 1 - e^{-\alpha_0}$;*

(b) *if $\pi_0 < 1$ and $A(\pi_0, \gamma) \leq \overline{A} < \infty$ for some $\overline{A}$ and sufficiently large $m$, then*

$$\Psi(\alpha^*_{cal}) \simeq \frac{\pi_0}{1 - \pi_0} \psi_*^{-1} \left[A(\pi_0, \gamma)\right]^{1-\xi_m} m^{-(1-\xi_m)B(\pi_0, \gamma)}, \quad as \; m \longrightarrow \infty.$$

*Proof.* See Appendix. □

There are two important consequences from this theorem. First, the level $\alpha_0$ can be chosen to bound ERR (and FDR) asymptotically in the least favorable situation $\pi_0 = 1$. In this case both ERR and FDR are equal to the family-wise type-I error probability. Note that $1 - e^{-\alpha_0}$ is also the limiting family-wise type-I error probability corresponding to the Bonferroni significance threshold $\alpha_0 m^{-1}$. In this regard the adaptive threshold $\alpha^*_{cal}$ is calibrated to the conservative Bonferroni threshold when $\pi_0 = 1$. If one wants to bound the error level at $\alpha_1$, then set $\alpha_0 = -\log(1 - \alpha_1)$. Of course $\alpha_0 \approx \alpha_1$ for small $\alpha_1$; for example, $\alpha_0 \approx 0.05129, 0.1054, 0.2231$ for $\alpha_1 = 0.05, 0.1, 0.2$ respectively.

Next, Part (b) demonstrates that if the "average power" of rejecting the false null hypotheses remains visible asymptotically in the sense that $\xi_m \leq \overline{\xi} < 1$ for some $\overline{\xi}$ and sufficiently large $m$, then the upper bound

$$\Psi(\alpha^*_{cal}) \simeq \frac{\pi_0}{1 - \pi_0} \psi_*^{-1} \left[A(\pi_0, \gamma)\right]^{1-\xi_m} m^{-(1-\overline{\xi})B(\pi_0, \gamma)} \longrightarrow 0;$$

therefore $ERR(\alpha^*_{cal})$ diminishes asymptotically. However, the convergence can be slow if the power is weak in the sense $\overline{\xi} \approx 1$ (hence $H_m(\cdot)$ is close to the $U(0,1)$ cdf in the left tail). Moreover, $\Psi$ can be considerably close to 1 in the unfavorable scenario $\pi_0 \approx 1$ and $\overline{\xi} \approx 1$. On the other hand, increase in the average power in the sense of decrease in $\overline{\xi}$ makes $\Psi$ (hence the ERR) diminishes faster asymptotically.

Note from (4.3) that as long as $\pi_0$ is bounded away from zero (i.e., there is always some null hypotheses remain true) and $\gamma$ is bounded, the quantity $A(\pi_0, \gamma)$ is bounded. Because the positive ERR does not involve the probability $\Pr(R > 0)$, part (b) holds for $pERR(\alpha^*_{cal})$ under *arbitrary* dependence among the P values (tests).

### 4.3. Data-driven adaptive significance threshold

Plugging $\widehat{\pi}_0$ and $\widehat{\gamma}$ generated by optimizing (3.3) into (4.2) produces a data-driven significance threshold:

$$(4.5) \qquad \widehat{\alpha}^*_{cal} := A(\widehat{\pi}_0, \widehat{\gamma}) m^{-B(\widehat{\pi}_0, \widehat{\gamma})}.$$

Now consider the ERR of the procedure $HT(\widehat{\alpha}^*_{cal})$ with $\widehat{\alpha}^*_{cal}$ as above. Define

$$ERR^* := \frac{E\left[V(\widehat{\alpha}^*_{cal})\right]}{E\left[R(\widehat{\alpha}^*_{cal})\right]} \Pr\left(R(\widehat{\alpha}^*_{cal}) > 0\right).$$

The interest here is the asymptotic magnitude of $ERR^*$ as $m \longrightarrow \infty$. A major difference here from Theorem 4.1 is that the threshold $\widehat{\alpha}^*_{cal}$ is random. A similar result can be established with some moment assumptions on $A(\widehat{\pi}_0, \widehat{\gamma})$, where $A(\cdot, \cdot)$ is defined in (4.3) and $\widehat{\pi}_0, \widehat{\gamma}$ are generated by optimizing (3.3). Toward this end, still assume that $\mathbf{P}_m$ is asymptotically stable, and let $\eta_m, \eta$, and $\xi_m$ be as in Definition 4.2. Let $\nu_m$ be the joint cdf of $[\widehat{\pi}_0, \widehat{\gamma}]$, and let

$$a_m := \int_{\mathbb{R}^2} A(s,t)^{\eta_m} d\nu_m(s,t)$$

$$a_{1m} := \int_{\mathbb{R}^2} A(s,t) d\nu_m(s,t)$$

$$a_{2m} := \int_{\mathbb{R}^2} A(s,t)^{\xi_m} d\nu_m(s,t).$$

All these moments exist as long as $\widehat{\pi}_0$ is bounded away from zero and $\widehat{\gamma}$ is bounded with probability one.

**Theorem 4.2.** *Suppose that* $\mathbf{P}_m$ *is asymptotically stable and has the positive orthant dependence property for sufficiently large* $m$. *Let* $\beta^*$, $\eta$, $\psi_*$, *and* $\xi_m$ *be as in Definition* 4.2. *If* $a_m$, $a_{1m}$ *and* $a_{2m}$ *all exist for sufficiently large* $m$, *then* $ERR^* \leq \Psi_m$ *and there exist* $\delta_m \in [\eta/3, \eta]$, $\varepsilon_m \in [1/3, 1]$, *and* $\varepsilon'_m \in [\xi_m/3, \xi_m]$ *such that as* $m \longrightarrow \infty$

$$\Psi_m \simeq \begin{cases} K(\beta^*, a_m, \delta_m), & \text{if } \pi_0 = 1, \text{ all } m \\ \frac{\pi_0}{1-\pi_0}\left(\frac{a_{1m}}{a_{2m}}\right)\psi_*^{-1}\frac{K(\beta^*, a_m, \delta_m)}{m^{(\varepsilon_m - \varepsilon'_m)}}, & \text{if } \pi_0 < 1, \text{ sufficiently large } m, \end{cases}$$

*where* $K(\beta^*, a_m, \delta_m) = 1 - \left(1 - \beta^* a_m m^{-\delta_m}\right)^m$

*Proof.* See Appendix. □

Although less specific than Theorem 4.1, this result still is instructive. First, if the "average power" sustains asymptotically in the sense that $\xi_m < 1/3$ so that $\varepsilon_m > \varepsilon'_m$ for sufficiently large $m$, or if $\lim_{m \to \infty} \xi_m = \overline{\xi} < 1/3$, then $ERR^*$ diminishes as $m \longrightarrow \infty$. The asymptotic behavior of $ERR^*$ in the case of $\overline{\xi} \geq 1/3$ is indefinite from this theorem, and obtaining a more detailed upper bound for $ERR^*$ in this case remains an open problem. Next, $ERR^*$ can be potentially high if $\pi_0 = 1$ always or $\pi_0 \approx 1$ and the average power is weak asymptotically. The reduced specificity in this result compared to Theorem 4.1 is due to the random variations in $A(\widehat{\pi}_0, \widehat{\gamma})$ and $B(\widehat{\pi}_0, \widehat{\gamma})$, which are now random variables instead of deterministic functions. Nonetheless Theorem 4.2 and its proof (see Appendix) do indicate that when $\pi_0 \approx 1$ and the average power is weak (i.e., $H_m(\cdot)$ is small), for sake of ERR (and FDR) reduction the variability in $A(\widehat{\pi}_0, \widehat{\gamma})$ and $B(\widehat{\pi}_0, \widehat{\gamma})$ should be reduced as much as possible in a way to make $\delta_m$ and $\varepsilon_m$ as close to 1 as possible. In practice one should make an effort to help this by setting $\widehat{\pi}_0$ and $\widehat{\gamma}$ to 1 when the smoothed empirical quantile function $\widehat{Q}_m$ is too close to the $U(0,1)$ quantile function. On the other hand, one would like to have a reasonable level of false negative errors when true alternative hypotheses do exist even if $\pi_0 \approx 1$; this can be helped by setting $\alpha_0$ at a reasonably liberal level. The simulation study (Section 5) indicates that $\alpha_0 = 0.22$ is a good choice in a wide variety of scenarios.

Finally note that, just like in Theorem 4.1, the bound when $\pi_0 < 1$ holds for the positive ERR $pERR^* := E[V(\widehat{\alpha}^*_{cal})]/E[R(\widehat{\alpha}^*_{cal})]$ under *arbitrary* dependence among the tests.

## 5. A Simulation study

To better understand and compare the performance and operating characteristics of $HT(\widehat{\alpha}^*_{cal})$, a simulation study is performed using models that mimic a gene signaling pathway to generate data, as proposed in [7]. Each simulation model is built from a network of 7 related "genes" (random variables), $X_0$, $X_1$, $X_2$, $X_3$, $X_4$, $X_{190}$, and $X_{221}$, as depicted in Figure 2, where $X_0$ is a latent variable. A number of other variables are linear functions of these random variables.

Ten models (scenarios) are simulated. In each model there are $m$ random variables, each observed in $K$ groups with $n_k$ independent observations in the $k$th group ($k = 1, \ldots K$). Let $\mu_{ik}$ be the mean of variable $i$ in group $k$. Then $m$ ANOVA hypotheses, one for each variable ($H_{0i}$: $\mu_{i1} = \cdots = \mu_{iK}$, $i = 1 \ldots, m$), are tested.

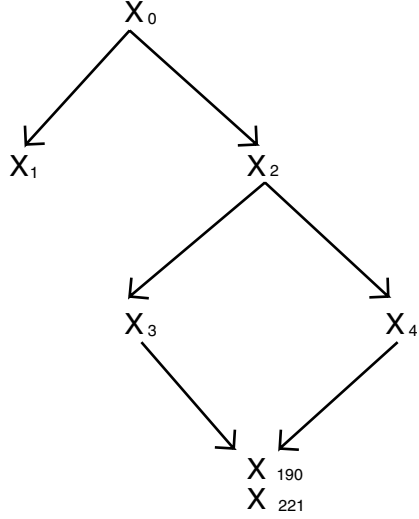

Fig 2. *A seven-variable framework to simulate differential gene expressions in a pathway.*

Table 2
*Relationships among $X_0, X_1, \ldots, X_4$, $X_{190}$ and $X_{221}$: $X_{ikj}$ denote the jth observation of the ith variable in group k; $N(0, \sigma^2)$ denotes normal random noise. The index j always runs through 1,2,3*

| |
|---|
| $X_{01j}$ i.i.d. $N(0, \sigma^2)$; $X_{0kj}$ i.i.d. $N(8, \sigma^2)$, $k = 2, 3, 4$ |
| $X_{1kj} = X_{0kj}/4 + N(0, 0.0784)$ ($X_1$ is highly correlated with $X_0$; $\sigma = 0.28$.) |
| $X_{2kj} = X_{0kj} + N(0, \sigma^2)$, $k = 1, 2$; $X_{23j} = X_{03j} + 6 + N(0, \sigma^2)$; $X_{24j} = X_{04j} + 14 + N(0, \sigma^2)$ |
| $X_{3kj} = X_{2kj} + N(0, \sigma^2)$, $k = 1, 2, 3, 4$ |
| $X_{4kj} = X_{2kj} + N(0, \sigma^2)$, $k = 1, 2$; $X_{43j} = X_{23j} - 6 + N(0, \sigma^2)$; $X_{44j} = X_{24j} - 8 + N(0, \sigma^2)$ |
| $X_{190,1j} = X_{31j} + 24 + N(0, \sigma^2)$; $X_{190,2j} = X_{32j} + X_{42j} + N(0, \sigma^2)$; |
| $\quad X_{190,3j} = X_{33j} - X_{43j} - 6 + N(0, \sigma^2)$; $X_{190,4j} = X_{34j} - 14 + N(0, \sigma^2)$ |
| $X_{221,kj} = X_{3kj} + 24 + N(0, \sigma^2)$, $k = 1, 2$; |
| $\quad X_{221,3j} = X_{33j} - X_{43j} + N(0, \sigma^2)$; $X_{221,4j} = X_{34j} + 2 + N(0, \sigma^2)$ |

Realizations are drawn from Normal distributions. For all ten models the number of groups $K = 4$ and the sample size $n_k = 3$, $k = 1, 2, 3, 4$. The usual one-way ANOVA $F$ test is used to calculate P values. Table 2 contains a detailed description of the joint distribution of $X_0, \ldots, X_4$, $X_{190}$ and $X_{221}$ in the ANOVA set up. The ten models comprised of different combinations of $m$, $\pi_0$, and the noise level $\sigma$ are detailed in Table 3, Appendix. The odd numbered models represent the high-noise (thus weak power) scenario and the even numbered models represent the low-noise (thus substantial power) scenario. In each model variables not mentioned in the table are i.i.d. $N(0, \sigma^2)$. Performance statistics under each model are calculated from 1,000 simulation runs.

First, the $\pi_0$ estimators by Benjamini and Hochberg [3], Storey *et al.* [31], and (3.1) are compared on several models. Root mean square error (MSE) and bias are plotted in Figure 3. In all cases the root MSE of the estimator (3.1) is either the smallest or comparable to the smallest. In the high noise case ($\sigma = 3$) Benjamini and Hochberg's estimator tends to be quite conservative (upward biased), especially for relatively low true $\pi_0$ (0.83 and 0.92, Models 1 and 3); whereas Storey's estimator is biased downward slightly in all cases. The proposed estimator (3.1) is biased in the conservative direction, but is less conservative than Benjamini and Hochberg's estimator. In the low noise case ($\sigma = 1$) the root MSE of all three estimators and the bias of the proposed and the Benjamini and Hochberg's estimators are reduced substantially while the small downward bias of Storey's bootstrap estimator remains. Overall the proposed estimator (3.1) outperforms the other two estimators in terms of MSE and bias.

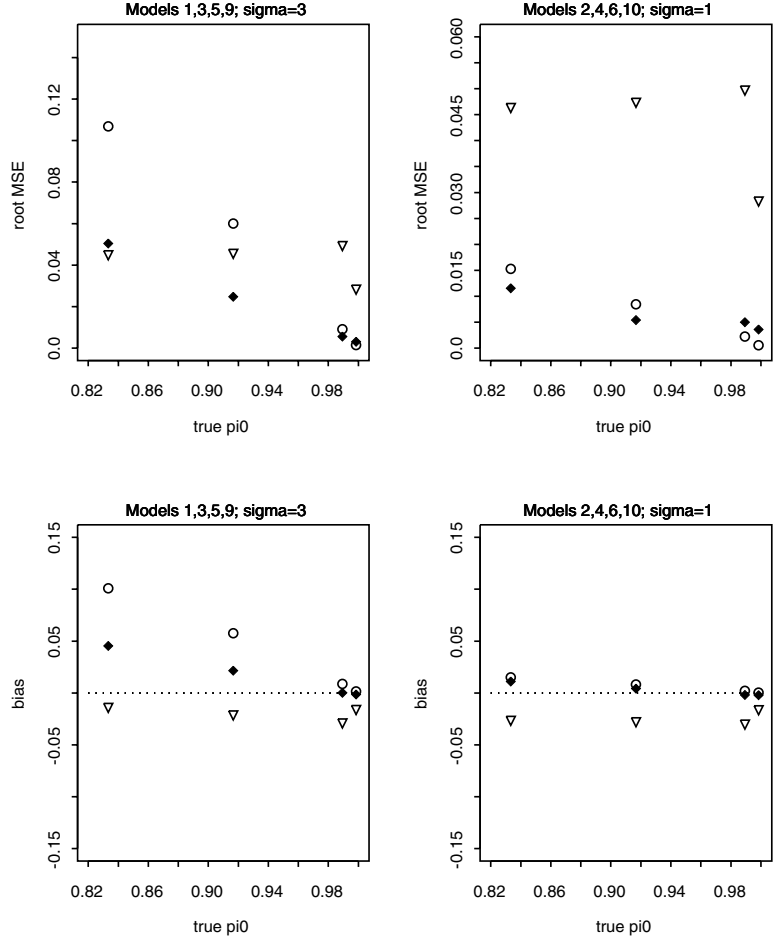

Fig 3. *Root MSE and bias of the $\pi_0$ estimators by Benjamini and Hochberg* [3] *(circle), Storey et al.* [31] *(triangle), and* (3.1) *(diamond)*

Next, operating characteristics of the adaptive FDR control ([3]) and $q$-value FDR control ([31]) at the 1%, 5%, 10%, 15%, 20%, 30%, 40%, 60%, and 70% levels, the criteria $API$ (i.e., the $HT(\widehat{\alpha}^*_{cal})$ procedure) and $I_p$ ([7]), are simulated and compared. The performance measures are the estimated FDR ($\widehat{FDR}$) and the estimated false nondiscovery proportion ($\widehat{FNDP}$) defined as follows. Let $m_1$ be the number of true alternative hypotheses according to the simulation model, let $R_l$ be the total number of rejections in simulation trial $l$, and let $S_l$ be the number of correct rejections. Define

$$\widehat{FDR} = \tfrac{1}{1000} \sum_{l=1}^{1000} I(R_l > 0)(R_l - S_l)/R_l$$
$$\widehat{FNDP} = \tfrac{1}{1000} \sum_{l=1}^{1000} (m_1 - S_l)/m_1,$$

where $I(\cdot)$ is the indicator function. These are the Monte Carlo estimators of the FDR and the $FNDP := E\left[m_1 - S\right]/m_1$ (cf. Table 1). In other words FNDP is the expected proportion of true alternative hypotheses not captured by the procedure. A measurement of the average power is $1 - FNDP$.

Following the discussions in Section 4, the parameter $\alpha_0$ required in the *API* procedure should be set at a reasonably liberal level. A few values of $\alpha_0$ were examined in a preliminary simulation study, which suggested that $\alpha_0 = 0.22$ is a level that worked well for the variety of scenarios covered by the ten models in Table 3, Appendix.

Results corresponding to $\alpha_0 = 0.22$ are reported here. The results are first summarized in Figure 4. In the high noise case ($\sigma = 3$, Models 1, 3, 5, 7, 9), compared to $I_p$, *API* incurs no or little increase in FNDP but substantially lower FDR when $\pi_0$ is high (Models 5, 7, 9), and keeps the same FDR level and a slightly reduced FNDP when $\pi_0$ is relatively low (Models 1, 3); thus *API* is more adaptive than $I_p$. As expected, it is difficult for all methods to have substantial power (low FNDP) in the high noise case, primarily due to the low power in each individual test to reject a false null hypothesis. For the FDR control procedures, no substantial number of false null hypotheses can be rejected unless the FDR control level is raised to a relatively high level of $\geq 30\%$, especially when $\pi_0$ is high.

In the low noise case ($\sigma = 1$, Models 2, 4, 6, 8, 10), *API* performs similarly to $I_p$, although it is slightly more liberal in terms of higher FDR and lower FNDP when $\pi_0$ is relatively low (Models 2, 4). Interestingly, when $\pi_0$ is high (Models 6, 8, 10), FDR control by $q$-value (Storey *et al.* [31]) is less powerful than the adaptive FDR procedure (Benjamini and Hochberg [3]) at low FDR control levels (1%, 5%, and 10%), in terms of elevated FNDP levels.

The methods are further compared by plotting $\widehat{FNDP}$ vs. $\widehat{FDR}$ for each model in Figure 5. The results demonstrate that in low-noise (model 2, 4, 6, 8, 10) and high-noise, high-$\pi_0$ (models 5, 7, 9) cases, the adaptive significance threshold determined from *API* gives very reasonable balance between the amounts of false positive and false negative errors, as indicated by the position of the diamond ($\widehat{FNDP}$ vs. $\widehat{FDR}$ of *API*) relative to the curves of the FDR-control procedures. It is noticeable that in the low noise cases the adaptive significance threshold corresponds well to the maximum FDR level for which there is no longer substantial gain in reducing FNDP by controlling the FDR at higher levels. There is some loss of efficiency for using *API* in high-noise, low-$\pi_0$ cases (model 1, 3) – its FNDP is higher than the control procedures at comparable FDR levels. This is a price to pay for not using a prespecified, fixed FDR control level.

The simulation results on *API* are very consistent with the theoretical results in Section 4. They indicate that *API* can provide a reasonable, data-adaptive significance threshold that balances the amounts of false positive and false negative errors: it is reasonably conservative in the high $\pi_0$ and high noise (hence low power) cases, and is reasonably liberal in the relatively low $\pi_0$ and low noise cases.

## 6. Concluding remarks

In this research an improved estimator of the null proportion and an adaptive significance threshold criterion *API* for massive multiple tests are developed and studied, following the introduction of a new measurement of the level of false positive errors, ERR, as an alternative to FDR for theoretical investigation. ERR allows for obtaining insights into the error behavior of *API* under more application-pertinent distributional assumptions that are widely satisfied by the data in many recent

70     *C. Cheng*

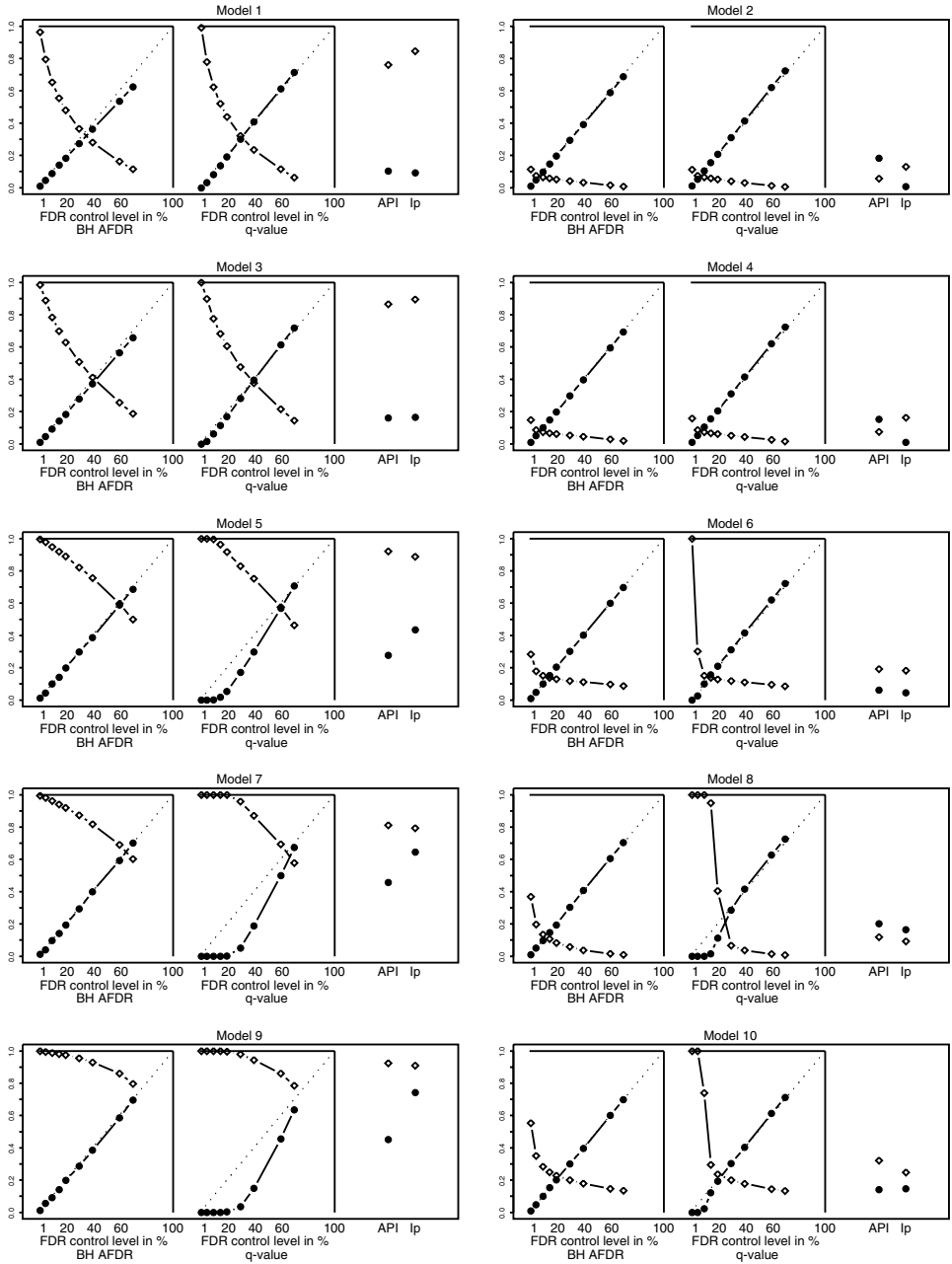

Fig 4. *Simulation results on the rejection criteria. Each panel corresponds to a model configuration. Panels in the left column correspond to the "high noise" case $\sigma = 3$, and panels in the right column correspond to the "low noise" case $\sigma = 1$. The performance statistics $\widehat{FDR}$ (bullet) and $\widehat{FNDP}$ (diamond) are plotted against each criteria. Each panel has three sections. The left section shows FDR control with the Benjamini & Hochberg [3] adaptive procedure (BH AFDR), and the middle section shows FDR control by q-value, all at the 1%, 5%, 10%, 15%, 20%, 30%, 40%, 60%, and 70% levels. The right section shows $\widehat{FDR}$ and $\widehat{FNDP}$ of API and $I_p$.*

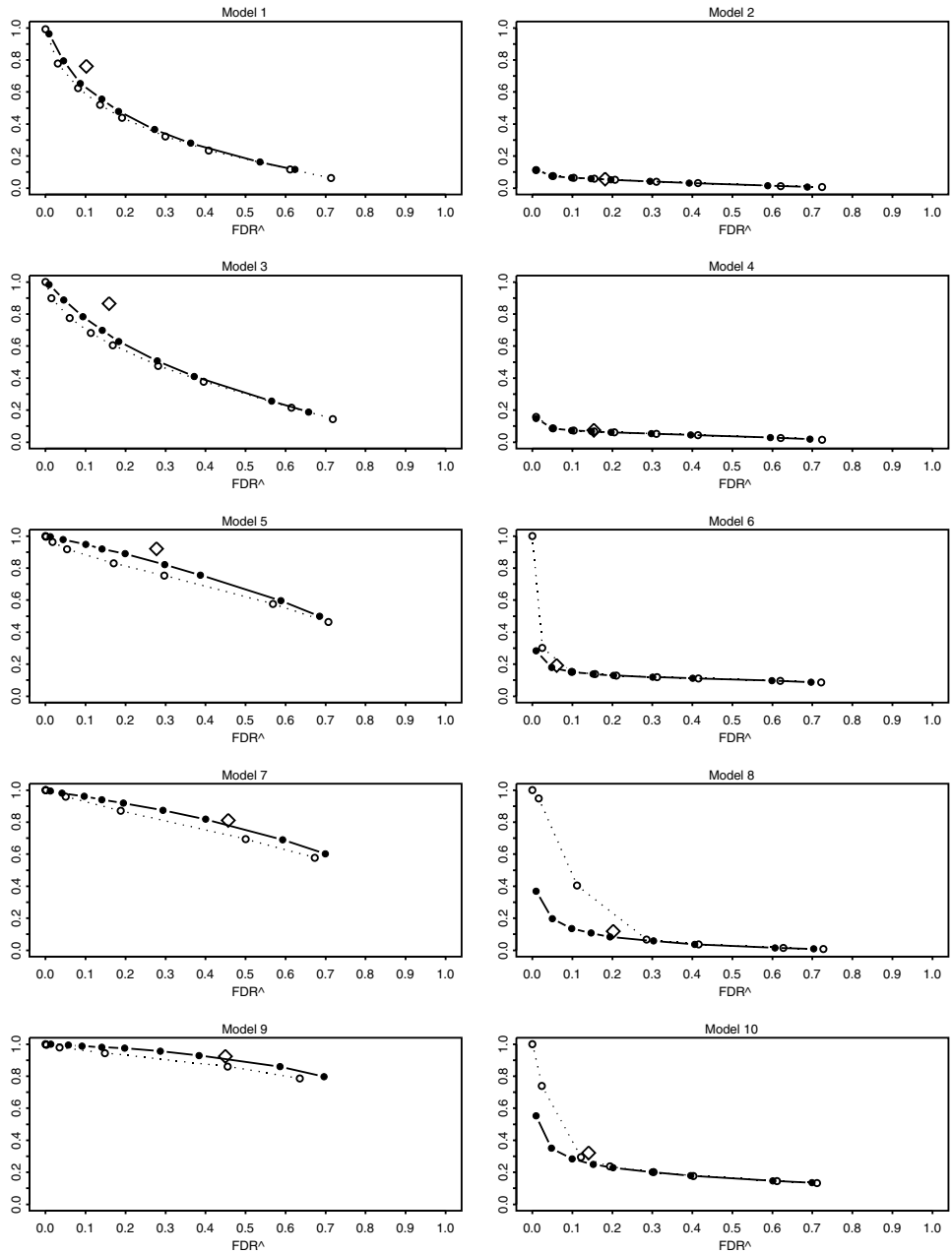

Fig 5. $\widehat{FNDP}$ vs. $\widehat{FDR}$ for Benjamini and Hochberg [3] adaptive FDR control (solid line and bullet) and q-value FDR control (dotted line and circle) when FDR control levels are set at 1%, 5%, 10%, 15%, 20%, 30%, 40%, 60%, and 70%. For each model $\widehat{FNDP}$ vs. $\widehat{FDR}$ of the adaptive API procedure occupies one point on the plot, indicated by a diamond.

applications. Under these assumptions, for the first time the asymptotic ERR level (and the FDR level under certain conditions) is explicitely related to the ensemble behavior of the P values described by the upper envelope cdf $\overline{F}_m$ and the "average power" $H_m$. Parallel to positive FDR, the concept of positive ERR is also useful. Asymptotic pERR properties of the proposed adaptive method can be established under arbitrary dependence among the tests. The theoretical understanding provides cautions and remedies to the application of $API$ in practice.

Under proper ergodicity conditions such as those used in [31, 14], FDR and ERR are equivalent for the hard-thresholding procedure (2.3); hence Theorems 4.1 and 4.2 hold for FDR as well.

The simulation study shows that the proposed estimator of the null proportion by quantile modeling is superior to the two popular estimators in terms of reduced MSE and bias. Not surprisingly, when there is little power to reject each individual false null hypothesis (hence little average power), FDR control and $API$ both incur high level of false negative errors in terms of FNDP. When there is a reasonable amount of power, $API$ can produce a reasonable balance between the false positives and false negatives, thereby complementing and extending the widely used FDR-control approach to massive multiple tests.

In exploratory type applications where it is desirable to provide "inference-guided discoveries", the role of $\alpha_0$ is to provide a protection in the situation where no true alternative hypothesis exits ($\pi_0 = 1$). On the other hand it is not advisable to choose the significance threshold too conservatively in such applications because the "discoveries" will be scrutinized in follow up investigations. Even if setting $\alpha_0 = 1$ the calibrated adaptive significance threshold is $m^{-1}$, giving the limiting ERR (or FDR, or family-wise type-I error probability) $1 - e^{-1} \approx 0.6321$ when $\pi_0 = 1$.

At least two open problems remain. First, although there has been empirical evidence from the simulation study that the $\pi_0$ estimator (3.1) outperforms the existing ones, there is lack of analytical understating of this estimator, in terms of MSE for example. Second, the bounds obtained in Theorems 4.1 and 4.2 are not sharp, a more detailed characterization of the upper bound of $ERR^*$ (Theorem 4.2) is desirable for further understanding of the asymptotic behavior of the adaptive procedure.

**Appendix**

*Proof of Theorem 4.1.* For (a), from (4.3) and (4.4), if $\pi_0 = 1$ for all $m$, then the first factor on the right-hand side of (4.4) is 1, and the second factor is now equal to

$$1 - (1 - \alpha_{cal}^*)^m = 1 - (1 - A(1,1)m^{-B(1,1)})^m = 1 - (1 - \alpha_0 m^{-1})^m \longrightarrow 1 - e^{-\alpha_0}$$

because $A(1,1) = \alpha_0$ and $B(1,1) = 1$. For (b), first

$$1 - (1 - \overline{F}_m(\alpha_{cal}^*))^m \simeq 1 - (1 - \beta_m \alpha_{cal}^{*\eta_m})^m \leq 1 - (1 - \beta^* \overline{A} \alpha_0 m^{-\eta B(\pi_0, \gamma)})^m := \varepsilon_m,$$

and $\varepsilon_m \simeq 1 - \exp\left(-\beta^* \overline{A} \alpha_0 m^{1-\eta B(\pi_0, \gamma)}\right) \longrightarrow 1$ because $B(\pi_0, \gamma) \leq (\gamma + 2)/(2\gamma + 1) < 1$ so that $\eta B(\pi_0, \gamma) < 1$. Next, let $\omega_m := \psi_m^{-1} \frac{[A(\pi_0, \gamma)]^{1-\xi_m}}{m^{(1-\xi_m)B(\pi_0, \gamma)}}$.

Then

$$\frac{\pi_0 \alpha_{cal}^*}{\pi_0 \alpha_{cal}^* + (1-\pi_0) H_m(\alpha_{cal}^*)} \simeq \frac{\pi_0 \omega_m}{1-\pi_0} \leq \frac{\pi_0}{1-\pi_0} \psi_*^{-1} \frac{[A(\pi_0, \gamma)]^{1-\xi_m}}{m^{(1-\xi_m)B(\pi_0, \gamma)}}$$

for sufficiently large $m$. Multiplying this upper bound and the limit of $\varepsilon_m$ gives (b). □

*Proof of Theorem 4.2.* First, for sufficiently large $m$,

$$\Pr\left(R(\widehat{\alpha}_{cal}^*) > 0\right) \leq 1 - \int_{\mathbb{R}^2} \left[\beta_m A(s,t)^{\eta_m} m^{-\eta_m B(s,t)}\right]^m d\nu_m(s,t)$$

$$\leq 1 - \left[1 - \beta^* \int_{\mathbb{R}^2} A(s,t)^{\eta_m} m^{-\eta B(s,t)} d\nu_m(s,t)\right]^m$$

Because $1/3 \leq B(\widehat{\pi}_0, \widehat{\gamma}) \leq 1$ with probability 1, so that $m^{-\eta} \leq m^{\eta B(\widehat{\pi}_0, \widehat{\gamma})} \leq m^{-\eta/3}$ with probability 1, by the mean value theorem of integration (Halmos [15], P.114), there exists some $\mathcal{E}_m \in [m^{-\eta}, m^{-\eta/3}]$ such that

$$\int_{\mathbb{R}^2} A(s,t)^{\eta_m} m^{-\eta B(s,t)} d\nu_m(s,t) = \mathcal{E}_m a_m,$$

and $\mathcal{E}_m$ can be written equivalently as $m^{-\delta_m}$ for some $\delta_m \in [\eta/3, \eta]$, giving, for sufficiently large $m$,

$$\Pr\left(R(\widehat{\alpha}_{cal}^*) > 0\right) \leq 1 - \left(1 - \beta^* a_m m^{-\delta_m}\right)^m.$$

This is the upper bound $\Psi_m$ of $ERR^*$ if $\pi_0 = 1$ for all $m$ because now $V(\widehat{\alpha}_{cal}^*) = R(\widehat{\alpha}_{cal}^*)$ with probability 1. Next,

$$E\left[V(\widehat{\alpha}_{cal}^*)\right] = E\left[E\left[V(\widehat{\alpha}_{cal}^*)\big|\widehat{\alpha}_{cal}^*\right]\right] = E\left[\pi_0 \widehat{\alpha}_{cal}^*\right] = \pi_0 \int_{\mathbb{R}^2} \frac{A(s,t)}{m^{B(s,t)}} d\nu_m(s,t).$$

Again by the mean value theorem of integration there exists $\varepsilon_m \in [1/3, 1]$ such that $E\left[V(\widehat{\alpha}_{cal}^*)\right] = \pi_0 a_{1m} m^{-\varepsilon_m}$. Similarly,

$$E\left[H_m(\widehat{\alpha}_{cal}^*)\right] \simeq \psi_m \int_{\mathbb{R}^2} A(s,t)^{\xi_m} m^{-\xi_m B(s,t)} d\nu_m(s,t) \geq \psi_* a_{2m} m^{-\varepsilon'_m}$$

for some $\varepsilon'_m \in [\xi_m/3, \xi_m]$. Finally, because

$$E\left[R(\widehat{\alpha}_{cal}^*)\right] = E\left[E\left[R(\widehat{\alpha}_{cal}^*)\big|\widehat{\alpha}_{cal}^*\right]\right] = \pi_0 E\left[V(\widehat{\alpha}_{cal}^*)\right] + (1-\pi_0) E\left[H_m(\widehat{\alpha}_{cal}^*)\right],$$

if $\pi_0 < 1$ for sufficiently large $m$, then

$$ERR^* \leq \left[1 - \left(1 - \beta^* a_m m^{-\delta_m}\right)^m\right] \frac{\pi_0 E[\widehat{\alpha}_{cal}^*]}{(1-\pi_0) E[H_m(\widehat{\alpha}_{cal}^*)]}$$

$$\leq \frac{\pi_0}{1-\pi_0} \left(\frac{a_{1m}}{a_{2m}}\right) \psi_*^{-1} m^{-(\varepsilon_m - \varepsilon'_m)} \left[1 - \left(1 - \beta^* a_m m^{-\delta_m}\right)^m\right]. \quad □$$

Table 3

*Ten models: Model configuration in terms of $(m, m_1, \sigma)$ and determination of true alternative hypotheses by $X_1, \ldots, X_4$, $X_{190}$ and $X_{221}$, where $m_1$ is the number of true alternative hypotheses; hence $\pi_0 = 1 - m_1/m$. $n_k = 3$ and $K = 4$ for all models. $N(0, \sigma^2)$ denotes normal random noise*

| Model | $m$ | $m_1$ | $\sigma$ | True $H_A$'s in addition to $X_1, \ldots X_4$, $X_{190}$, $X_{221}$ |
|---|---|---|---|---|
| 1 | 3000 | 500 | 3 | $X_i = X_1 + N(0, \sigma^2)$, $i = 5, \ldots, 16$ |
| | | | | $X_i = -X_1 + N(0, \sigma^2)$, $i = 17, \ldots, 25$ |
| | | | | $X_i = X_2 + N(0, \sigma^2)$, $i = 26, \ldots, 60$ |
| | | | | $X_i = -X_2 + N(0, \sigma^2)$, $i = 61, \ldots, 70$ |
| | | | | $X_i = X_3 + N(0, \sigma^2)$, $i = 71, \ldots, 100$ |
| | | | | $X_i = -X_3 + N(0, \sigma^2)$, $i = 101, \ldots, 110$ |
| | | | | $X_i = X_4 + N(0, \sigma^2)$, $i = 111, \ldots, 150$ |
| | | | | $X_i = -X_4 + N(0, \sigma^2)$, $i = 151, \ldots, 189$ |
| | | | | $X_i = X_{190} + N(0, \sigma^2)$, $i = 191, \ldots, 210$ |
| | | | | $X_i = -X_{190} + N(0, \sigma^2)$, $i = 211, \ldots, 220$ |
| | | | | $X_i = X_{221} + N(0, \sigma^2)$, $i = 222, \ldots, 250$ |
| | | | | $X_i = 2X_{i-250} + N(0, \sigma^2)$, $i = 251, \ldots, 500$ |
| 2 | 3000 | 500 | 1 | the same as Model 1 |
| 3 | 3000 | 250 | 3 | the same as Model 1 except only the first 250 are true $H_A$'s |
| 4 | 3000 | 250 | 1 | the same as Model 3 |
| 5 | 3000 | 32 | 3 | $X_i = X_1 + N(0, \sigma^2)$, $i = 5, \ldots, 8$ |
| | | | | $X_i = X_2 + N(0, \sigma^2)$, $i = 9, \ldots, 12$ |
| | | | | $X_i = X_3 + N(0, \sigma^2)$, $i = 13, \ldots, 16$ |
| | | | | $X_i = X_4 + N(0, \sigma^2)$, $i = 17, \ldots, 20$ |
| | | | | $X_i = X_{190} + N(0, \sigma^2)$, $i = 191, \ldots, 195$ |
| | | | | $X_i = X_{221} + N(0, \sigma^2)$, $i = 222, \ldots, 226$ |
| 6 | 3000 | 32 | 1 | the same as Model 5 |
| 7 | 3000 | 6 | 3 | none, except $X_1, \ldots X_4$, $X_{190}$, $X_{221}$ |
| 8 | 3000 | 6 | 1 | the same as Model 7 |
| 9 | 10000 | 15 | 3 | $X_i = X_1 + N(0, \sigma^2)$, $i = 5, 6$ |
| | | | | $X_i = X_2 + N(0, \sigma^2)$, $i = 7, 8$ |
| | | | | $X_i = X_3 + N(0, \sigma^2)$, $i = 9, 10$ |
| | | | | $X_i = X_4 + N(0, \sigma^2)$, $i = 11, 12$ |
| | | | | $X_{191} = X_{190} + N(0, \sigma^2)$ |
| 10 | 10000 | 15 | 1 | the same as Model 9 |

**Acknowledgments.** I am grateful to Dr. Stan Pounds, two referees, and Professor Javier Rojo for their comments and suggestions that substantially improved this paper.

**References**

[1] ABRAMOVICH, F., BENJAMINI, Y., DONOHO, D. AND JOHNSTONE, I. (2000). Adapting to unknown sparsity by controlling the false discover rate. Technical Report 2000-19, Department of Statistics, Stanford University, Stanford, CA.

[2] ALLISON, D. B., GADBURY, G. L., HEO, M. FERNANDEZ, J. R., LEE, C-K, PROLLA, T. A. AND WEINDRUCH, R. (2002). A mixture model approach for the analysis of microarray gene expression data. *Comput. Statist. Data Anal.* **39**, 1–20. MR1895555

[3] BENJAMINI, Y. AND HOCHBERG, Y. (2000). On the adaptive control of the false discovery rate in multiple testing with independent statistics. *Journal of Educational and Behavioral Statistics* **25**, 60–83.

[4] BENJAMINI, Y. AND HOCHBERG, Y. (1995). Controlling the false discovery rate: a practical and powerful approach to multiple testing. *J. R. Stat. Soc. Ser. B Stat. Methodol.* **57**, 289–300. MR1325392

[5] BENJAMINI, Y., KRIEGER, A. M. AND YEKUTIELI, D. (2005). Adaptive linear step-up procedures that control the false discovery rate. Research Paper 01-03, Dept. of Statistics and Operations Research, Tel Aviv University.

[6] BICKEL, D. R. (2004). Error-rate and decision-theoretic methods of multiple testing: which genes have high objective probabilities of differential expression? *Statistical Applications in Genetics and Molecular Biology* **3**, Article 8. URL //www.bepress.com/sagmb/vol3/iss1/art8. MR2101458

[7] CHENG, C., POUNDS, S., BOYETT, J. M., PEI, D., KUO, M-L., ROUSSEL, M. F. (2004). Statistical significance threshold criteria for analysis of microarray gene expression data. *Statistical Applications in Genetics and Molecular Biology* **3**, Article 36. URL //www.bepress.com/sagmb/vol3/iss1/art36. MR2101483

[8] DE BOOR (1987). *A Practical Guide to Splines.* Springer, New York.

[9] DUDOIT, S., VAN DER LAAN, M., POLLARD, K. S. (2004). Multiple Testing. Part I. Single-step procedures for control of general Type I error rates. *Statistical Applications in Genetics and Molecular Biology* **3**, Article 13. URL //www.bepress.com/sagmb/vol3/iss1/art13. MR2101462

[10] EFRON, B. (2004). Large-scale simultaneous hypothesis testing: the choice of a null hypothesis. *J. Amer. Statist. Assoc.* **99**, 96–104. MR2054289

[11] EFRON, B., TIBSHIRANI, R., STOREY, J. D. AND TUSHER, V. (2001). Empirical Bayes analysis of a microarray experiment. *J. Amer. Statist. Assoc.* **96**, 1151–1160. MR1946571

[12] FINNER, H. AND ROBERTS, M. (2002). Multiple hypotheses testing and expected number of type I errors. *Ann. Statist.* **30**, 220–238. MR1892662

[13] GENOVESE, C. AND WASSERMAN, L. (2002). Operating characteristics and extensions of the false discovery rate procedure. *J. R. Stat. Soc. Ser. B Stat. Methodol.* **64**, 499–517. MR1924303

[14] GENOVESE, C. AND WASSERMAN, L. (2004). A stochastic process approach to false discovery rates. *Ann. Statist.* **32**, 1035–1061. MR2065197

[15] HALMOS, P. R. (1974). *Measure Theory.* Springer, New York. MR0453532

[16] HARDY, G., LITTLEWOOD, J. E. AND PÓLYA, G. (1952). *Inequalities.* Cambridge University Press, Cambridge, UK. MR0046395

[17] ISHAWARAN, H. AND RAO, S. (2003). Detecting differentially genes in microarrays using Baysian model selection. *J. Amer. Statist. Assoc.* **98**, 438–455. MR1995720

[18] KUO, M.-L., DUNCAVICH, E., CHENG, C., PEI, D., SHERR, C. J. AND ROUSSEL M. F. (2003). Arf induces p53-dependent and in-dependent antiproliferative genes. *Cancer Research* **1**, 1046–1053.

[19] LANGAAS, M., FERKINGSTADY, E. AND LINDQVIST, B. H. (2005). Estimating the proportion of true null hypotheses, with application to DNA microarray data. *J. R. Stat. Soc. Ser. B Stat. Methodol.* **67**, 555–572. MR2168204

[20] LEHMANN, E. L. (1966). Some concepts of dependence. *Ann. Math. Statist* **37**, 1137–1153. MR0202228

[21] MOSIG, M. O., LIPKIN, E., GALINA, K., TCHOURZYNA, E., SOLLER, M. AND FRIEDMANN, A. (2001). A whole genome scan for quantitative trait loci affecting milk protein percentage in Israeli-Holstein cattle, by means of selective milk DNA pooling in a daughter design, using an adjusted false discovery rate criterion. *Genetics* **157**, 1683–1698.

[22] NETTLETON, D. AND HWANG, G. (2003). Estimating the number of false null hypotheses when conducting many tests. Technical Report 2003-09, Department of Statistics, Iowa State University,

```
http://www.stat.iastate.edu/preprint/articles/2003-09.pdf
```
[23] NEWTON, M. A., NOUEIRY, A., SARKAR, D. AND AHLQUIST, P. (2004). Detecting differential gene expression with a semiparametric hierarchical mixture method, *Biostatistics* **5**, 155–176.

[24] POUNDS, S. AND MORRIS, S. (2003). Estimating the occurrence of false positives and false negatives in microarray studies by approximating and partitioning the empirical distribution of $p$-values. *Bioinformatics* **19**, 1236–1242.

[25] POUNDS, S. AND CHENG, C. (2004). Improving false discovery rate estimation. *Bioinformatics* **20**, 1737–1745.

[26] REINER, A., YEKUTIELI, D. AND BENJAMINI, Y. (2003). Identifying differentially expressed genes using false discovery rate controlling procedures. *Bioinformatics* **19**, 368–375.

[27] SCHWEDER, T. AND SPJØTVOLL, E. (1982). Plots of $P$-values to evaluate many tests simultaneously. *Biometrika* **69** 493-502.

[28] SMYTH, G. K. (2004). Linear models and empirical Bayes methods for assessing differential expression in microarray experiments. *Statistical Applications in Genetics and Molecular Biology* **3**, Article 3. URL: //www.bepress.com/sagmb/vol3/iss1/art3. MR2101454

[29] STOREY, J. D. (2002). A direct approach to false discovery rates. *J. R. Stat. Soc. Ser. B Stat. Methodol.* **64**, 479–498. MR1924302

[30] STOREY, J. D. (2003). The positive false discovery rate: a Baysian interpretation and the $q$-value. *Ann. Statis.* **31**, 2103–2035. MR2036398

[31] STOREY, J. D., TAYLOR, J. AND SIEGMUND, D. (2003). Strong control, conservative point estimation and simultaneous conservative consistency of false discovery rates: a unified approach. *J. R. Stat. Soc. Ser. B Stat. Methodol.* **66**, 187–205. MR2035766

[32] STOREY, J. D. AND TIBSHIRANI, R. (2003). SAM thresholding and false discovery rates for detecting differential gene expression in DNA microarrays. In *The Analysis of Gene Expression Data* (Parmigiani, G. et al., eds.). Springer, New York. MR2001400

[33] STOREY, J. D. AND TIBSHIRANI, R. (2003). Statistical significance for genome-wide studies. *Proc. Natl. Acad. Sci. USA* **100**, 9440–9445. MR1994856

[34] TSAI, C-A., HSUEH, H-M. AND CHEN, J. J. (2003). Estimation of false discovery rates in multiple testing: Application to gene microarray data. *Biometrics* **59**, 1071–1081. MR2025132

[35] TUSHER, V. G., TIBSHIRANI, R. AND CHU, G. (2001). Significance analysis of microarrays applied to ionizing radiation response. *Proc. Natl. Acad. Sci. USA* **98**, 5116–5121.

[36] VAN DER LAAN, M., DUDOIT, S. AND POLLARD, K. S. (2004a). Multiple Testing. Part II. Step-down procedures for control of the family-wise error rate. *Statistical Applications in Genetics and Molecular Biology* **3**, Article 14. URL: //www.bepress.com/sagmb/vol3/iss1/art14. MR2101463

[37] VAN DER LAAN, M., DUDOIT, S. AND POLLARD, K. S. (2004b). Augmentation procedures for control of the generalized family-wise error rate and tail probabilities for the proportion of false positives. *Statistical Applications in Genetics and Molecular Biology* **3**, Article 15. URL: //www.bepress.com/sagmb/vol3/iss1/art15. MR2101464